\documentclass{article}
\usepackage{amsmath}
\usepackage{amsfonts}

\newtheorem{theorem}{Theorem}

\newtheorem{corollary}[theorem]{Corollary}

\newtheorem{definition}[theorem]{Definition}

\newtheorem{lemma}[theorem]{Lemma}

\newenvironment{proof}[1][Proof]{\noindent\textbf{#1} }{\ \rule{0.5em}{0.5em}}
\input{tcilatex}
\begin{document}

\begin{center}
\textbf{Euler Coefficients and Restricted Dyck Paths}

Heinrich Niederhausen and Shaun Sullivan,

Department of Mathematics, Florida Atlantic University, Boca Raton\bigskip 
\end{center}

We consider the problem of enumerating Dyck paths staying weakly above the $x
$-axis with a limit to the number of consecutive $\searrow $ steps, or a
limit to the number of consecutive $\nearrow $ steps. We use Finite Operator
Calculus to obtain formulas for the number of all such paths reaching a
given point in the first quadrant. All our results are based on the Eulerian
coefficients.

\section{Introduction}

One of the most recent papers on patterns occurring $k$ times in Dyck paths
was written by A. Sapounakis, I. Tasoulas, P. Tsikouras, Counting strings in
Dyck paths, 2007, to appear in \emph{Discrete Mathematics} \cite{Sapounakis}.%
\textit{\ }The authors find generating functions for all $16$ patterns
generated by combinations of four up ($\nearrow $) and down ($\searrow $)
steps. A Dyck path starts at $\left( 0,0\right) $, takes only up and down
steps, and ends at $\left( 2n,0\right) $, staying weakly above the $x$-axis.
Returning to the $x$-axis at the end of the path has the advantage that
every path containing the pattern $uduu$, say, $k$ times, will contain the
reversed pattern $ddud$ also $k$ times when read backwards. This reduces
significantly the number of patterns under consideration. Dyck paths
containing $k$ strings of length $3$ were discussed by E. Deutsch in \cite%
{Deutsch}.

In this paper we consider only the patterns $u^{r}$ and $d^{r}$, for all
integers $r>2$, and we will investigate only the case $k=0$, which means
pattern avoidance. It has been shown in \cite{Sapounakis} that the
generating function $f(t)$ for avoiding $u^{r}$(or $d^{r}$) satisfies the
equation $f\left( t\right) =1+\sum_{i=1}^{r-1}t^{i}f\left( t\right) ^{i}=%
\frac{1-t-t^{r}f\left( t\right) ^{r}}{1-2t}$. However, we will allow the
Dyck paths to end at $\left( n,m\right) $, $m\geq 0$, which removes the
above mentioned symmetry, as shown in the following two tables.

\begin{equation*}
\begin{tabular}{l}
$%
\begin{tabular}{l||llllllllllllll}
$m$ &  &  &  &  &  &  &  &  &  &  &  &  &  &  \\ 
$7$ &  &  &  &  &  &  &  &  &  &  &  & 1 &  & 10 \\ 
$6$ &  &  &  &  &  &  &  &  &  &  & 3 &  & 19 &  \\ 
$5$ &  &  &  &  &  &  &  & 1 &  & 6 &  & 28 &  & 112 \\ 
$4$ &  &  &  &  &  &  & 2 &  & 9 &  & 33 &  & 116 &  \\ 
$3$ &  &  &  & 1 &  & 3 &  & 10 &  & 32 &  & 101 &  & 321 \\ 
$2$ &  &  & 1 &  & 3 &  & 8 &  & 23 &  & 68 &  & 205 &  \\ 
$1$ &  & 1 &  & 2 &  & 5 &  & 13 &  & 36 &  & 104 &  & 309 \\ 
$0$ & 1 &  & 1 &  & 2 &  & 5 &  & 13 &  & 36 &  & 104 &  \\ \hline\hline
$n:$ & $0$ & $1$ & $2$ & $3$ & $4$ & $5$ & $6$ & $7$ & $8$ & $9$ & $10$ & $%
11 $ & $12$ & $13$%
\end{tabular}%
$ \\ 
\multicolumn{1}{c}{The number of Dyck paths avoiding $uuuu$}%
\end{tabular}%
\end{equation*}

\begin{equation*}
\begin{tabular}{l}
$%
\begin{tabular}{l||llllllllllllll}
$m$ &  &  &  &  &  &  &  &  &  &  &  &  &  &  \\ 
$7$ &  &  &  &  &  &  &  & 1 &  & 8 &  & 44 &  & 208 \\ 
$6$ &  &  &  &  &  &  & 1 &  & 7 &  & 35 &  & 154 &  \\ 
$5$ &  &  &  &  &  & 1 &  & 6 &  & 27 &  & 110 &  & 423 \\ 
$4$ &  &  &  &  & 1 &  & 5 &  & 20 &  & 75 &  & 270 &  \\ 
$3$ &  &  &  & 1 &  & 4 &  & 14 &  & 48 &  & 161 &  & 536 \\ 
$2$ &  &  & 1 &  & 3 &  & 9 &  & 28 &  & 87 &  & 273 &  \\ 
$1$ &  & 1 &  & 2 &  & 5 &  & 14 &  & 40 &  & 118 &  & 357 \\ 
$0$ & 1 &  & 1 &  & 2 &  & 5 &  & 13 &  & 36 &  & 104 &  \\ \hline\hline
$n:$ & $0$ & $1$ & $2$ & $3$ & $4$ & $5$ & $6$ & $7$ & $8$ & $9$ & $10$ & $%
11 $ & $12$ & $13$%
\end{tabular}%
$ \\ 
\multicolumn{1}{c}{The number of Dyck paths avoiding $dddd$}%
\end{tabular}%
\end{equation*}%
The two tables indicate the differences between the two problems, both
starting out from equal counts on the $x$-axis ($m=0$). Because only points $%
\left( n,m\right) $ with $n+m=0\ \func{mod}2$ can be reached by a Dyck path,
we consider the lattice points $\left( 2n+b,2m+b\right) $, for $b=0,1$. We
first show that the number of Dyck paths to $\left( 2n+b,2m+b\right) $
avoiding $d^{r}$ equals 
\begin{equation*}
Dyck\left( 2n+b,2m+b;d^{r}\right) =\dfrac{2m+b+1}{n+m+b+1}\dbinom{n+m+b+1}{%
n-m}_{r},
\end{equation*}%
where the Euler coefficient \cite{Euler} is denoted by 
\begin{equation*}
\dbinom{n+m+b+1}{n-m}_{r}=\sum\limits_{i=0}^{\left\lfloor \left( n-m\right)
/r\right\rfloor }(-1)^{i}\dbinom{n+m+b+1}{i}\dbinom{2n+b-ri}{n-m-ri}
\end{equation*}%
(see Definition \ref{DefEuler} and expansion (\ref{(EulerBinomial)})). More
about Euler coefficients can be found in Section \ref{Euler}. For given $m$,
the number of Dyck paths \newline
$Dyck\left( n+m,m-n;d^{r}\right) $ to $\left( n+m,m-n\right) $ avoiding $%
d^{r}$ has the generating function (over $n$) $\left( 1-t^{r}\right)
^{m}\left( 1-t\right) ^{-m-2}\left( rt^{r}\left( 1-t\right) -\left(
1-t^{r}\right) \left( 2t-1\right) \right) $, as shown in (\ref{(GenFuncDown)}%
). Note that the coefficient of $t^{m}$ in this generating function equals
the original $Dyck\left( 2m,0;d^{r}\right) $.

Next we show that the number of Dyck paths to $\left( 2n+b,2m+b\right) $
avoiding $u^{r}$ equals 
\begin{eqnarray*}
&&D\left( 2m+b,2m+b;u^{4}\right) \\
&=&\sum_{i=0}^{2m+b-1}\frac{1}{n+m+b+1-i}\dbinom{i-2m-b}{i}_{r}\dbinom{%
n+m+b+1-i}{n+m+b-i}_{r},
\end{eqnarray*}
except for the original Dyck path counts to $\left( 2n,0\right) $, which
either must be gotten from those to $\left( 2n-1,1\right) $, or from the
Dyck paths to $\left( 2n,0\right) $ avoiding $d^{r}$. The case $r=4$ seems
to be very special. We conjecture in Section \ref{Conjecture} that in this
case the generating function for the Dyck paths equals%
\begin{eqnarray*}
&&\sum_{n\geq 0}Dyck\left( 4m-n-1,2m-n+1;u^{4}\right) \\
&=&\left( 3+t-\sqrt{\left( 1+t\right) ^{2}+4t^{3}}\right) \left( \frac{%
1-t^{4}}{1-t}\right) ^{m}/2,
\end{eqnarray*}%
hence $Dyck\left( 2m,0;u^{4}\right) =\left[ t^{2m}\right] \left( 3+t-\sqrt{%
\left( 1+t\right) ^{2}+4t^{3}}\right) \left( \frac{1-t^{4}}{1-t}\right)
^{m}/2$.

Throughout the following sections we will discuss ballot paths (weakly above 
$y=x$), with steps $\uparrow $ and $\rightarrow $, instead of Dyck paths.
The transformations $D\left( n,m\right) =Dyck\left( n+m,m-n\right) $ and $%
Dyck\left( 2n+b,2m+b\right) =D\left( n-m,n+m+b\right) $, with $D\left(
n,m\right) $ counting ballot path to $\left( n,m\right) $, go back and forth
between the two equivalent setups. Of course, the pattern $u^{r}$ becomes
the pattern $\uparrow ^{r}$, or $N^{r}$, and $d^{r}$ becomes $\rightarrow
^{r}$, or $E^{r}$.

\section{Ballot paths without the pattern $\rightarrow ^{r}$}

\begin{definition}
$s_{n}(m;r)=s_{n}(m)$ is the number of $\left\{ \uparrow ,\rightarrow
\right\} $ paths staying weakly above the diagonal $y=x$ from $(0,0)$ to $%
(n,m)\in \mathbb{Z}^{2}$ avoiding a sequence of $r>0$ consecutive $%
\rightarrow $ steps. We get $s_{0}\left( m\right) =1$ for all $m\geq 0$. We
set $s_{n}(m)=0$ if $n<0$ or if $m+1=n>0$.
\end{definition}

\begin{lemma}
The following recurrence holds for all $m\geq n>0$: 
\begin{equation}
s_{n}(m)=s_{n-1}(m)+s_{n}(m-1)-s_{n-r}(m-1).  \label{(MainRecurrence)}
\end{equation}
\end{lemma}

\begin{proof}[Proof:]
The number of paths reaching $(n,m)$ is obtained by adding the number of
paths reaching $(n-1,m)$ and $(n,m-1)$, but subtracting paths that would
have exactly $r$ $\rightarrow $ steps. Those forbidden steps occur
necessarily at the end of the path, so they are preceded by an up step, and
must come from $(n-r,m-1)$.
\end{proof}

We now extend $s_{n}(m)$ to all integers $m$ by first setting $s_{0}(m)=1$
and using (\ref{(MainRecurrence)}) to define the remaining $s_{n}(m)$ for $%
m<n-1$.

\begin{lemma}
$(s_n)$ is a polynomial sequence with $\deg{s_n}=n$.
\end{lemma}

\begin{proof}[Proof:]
We proceed by induction on $n$. Clearly, $\deg \left( s_{0}\right) =0$.
Suppose $s_{k}(m)$ is a polynomial of degree $k$ for $0\leq k\leq l$. Then $%
s_{l+1}(m)-s_{l+1}(m-1)=s_{l}(m)-s_{l-r+1}(m-1)$, which implies the first
difference of $s_{l+1}(m)$ is a polynomial of degree $l$. Thus, $s_{l+1}(m)$
is a polynomial of degree $l+1$.
\end{proof}

By interpolation we can define $\left( s_{n}\right) $ on all real numbers.

\begin{definition}
\label{DefEuler}The Eulerian Coefficient is defined as

\begin{eqnarray*}
\dbinom{x}{n}_{r} &=&[t^{n}](1+t+\cdots +t^{r-1})^{x} \\
&=&\sum\limits_{i=0}^{\left\lfloor n/r\right\rfloor }(-1)^{i}\dbinom{x}{i}%
\dbinom{x+n-ri-1}{n-ri}
\end{eqnarray*}%
(see (\ref{(EulerBinomial)})). Note that for $r=2$ the Euler coefficient
equals the binomial coefficient $\dbinom{x}{n}$.
\end{definition}

The following table shows the polynomial extension of $s_{n}(m)$. The number
of $\left\{ \uparrow ,\rightarrow \right\} $ paths to $(n,m)$ avoiding a
sequence of 4 $\rightarrow $ steps appear above the $y=x$ diagonal. The
numbers on the diagonal $\left( n,n\right) $, $1,1,2,5,13,36,\dots $, are
the number of Dyck paths to $\left( 2n,0\right) $. Of course, $s_{n}\left(
n\right) \leq C_{n}$, the $n$-th Catalan number.

\begin{equation*}
\begin{tabular}{c}
$%
\begin{tabular}{r||rrrrrrrrr}
$m$ & 1 & 7 & 27 & 75 & 161 & 273 & 357 & 309 & 0 \\ 
6 & 1 & 6 & 20 & 48 & 87 & 118 & 104 & 0 & -222 \\ 
5 & 1 & 5 & 14 & 28 & 40 & 36 & 0 & -76 & -182 \\ 
4 & 1 & 4 & 9 & 14 & 13 & 0 & -27 & -62 & -93 \\ 
3 & 1 & 3 & 5 & 5 & 0 & -10 & -22 & -30 & -31 \\ 
2 & 1 & 2 & 2 & 0 & -4 & -8 & -10 & -8 & -5 \\ 
1 & 1 & 1 & 0 & -2 & -3 & -3 & -2 & 0 & 0 \\ 
0 & 1 & 0 & -1 & -2 & 0 & 0 & 0 & 0 & 0 \\ 
-1 & 1 & -1 & -1 & -1 & 3 & -1 & -1 & -1 & 3 \\ \hline\hline
$n:$ & 0 & 1 & 2 & 3 & 4 & 5 & 6 & 7 & 8%
\end{tabular}%
$ \\ 
The path counts $s_{n}\left( m\right) $ and their polynomial extension ($r=4$%
)%
\end{tabular}%
\end{equation*}

\begin{theorem}
\label{ThDownSteps}%
\begin{equation*}
s_{n}(x)=\dfrac{x-n+1}{x+1}\dbinom{x+1}{n}_{r}=\dfrac{x-n+1}{x+1}%
\sum\limits_{i=0}^{\left\lfloor n/r\right\rfloor }(-1)^{i}\dbinom{x+1}{i}%
\dbinom{x+n-ri}{n-ri}
\end{equation*}
\end{theorem}

\begin{proof}[Proof:]
We saw that $\left( s_{n}\left( x\right) \right) $ is a basis for the vector
space of polynomials. Using operators on polynomials, we can write the
recurrence relation as%
\begin{equation}
1-E^{-1}=B-B^{r}E^{-1}  \label{(OperatorRec1)}
\end{equation}%
where $B$ and $E^{a}$ are defined by linear extension of $%
Bs_{n}(x)=s_{n-1}(x)$ and $E^{a}s_{n}(x)=s_{n}(x+a)$, the shift by $a$. The
operators $\nabla =1-E^{-1}$ and $E^{-1}$ both have power series expansions
in $D$, the derivative operator. Hence $B$ must have such an expansion too,
and therefore commutes with $\nabla $ and $E^{a}$. The power series for $B$
must be of order $1$, because $B$ reduces degrees by $1$. Such linear
operators are called \textit{delta operators}. The basic sequence $\left(
b_{n}\left( x\right) \right) _{n\geq 0}$ of a delta operator $B$ is a
sequence of polynomials such that $\deg b_{n}=n$, $Bb_{n}\left( x\right)
=b_{n-1}\left( x\right) $ (like the \textit{Sheffer sequence} $s_{n}\left(
x\right) $ for $B$), and initial conditions $b_{n}\left( 0\right) =\delta
_{0,n}$for all $n\in \mathbb{N}_{0}$. In our special case, the basic
sequence is easily determined. Solving for $E^{1}$ in (\ref{(OperatorRec1)})
shows that 
\begin{equation*}
E^{1}=\sum\limits_{i=0}^{r-1}B^{i}.
\end{equation*}%
Finite Operator Calculus tells us that if $E^{1}=1+\sigma (B)$, where $%
\sigma (t)$ is a power series of order 1 \cite[(2.5)]{Nied03}, then the
basic sequence $b_{n}(x)$ of $B$ has the generating function

\begin{equation*}
\sum\limits_{n\geq 0}b_{n}(x)t^{n}=(1+\sigma (t))^{x}.
\end{equation*}%
Thus, in our case $b_{n}(x)=\left[ t^{n}\right] (1+t+t^{2}+\cdots
+t^{r-1})^{x}=\dbinom{x}{n}_{r}$. Since the Sheffer sequence $(s_{n})$ has
initial values $s_{n}(n-1)=\delta _{n,0}$, using Abelization \cite{Nied03}
gives us

\begin{equation}
s_{n}(x)=\dfrac{x-n+1}{x+1}b_{n}(x+1)=\dfrac{x-n+1}{x+1}\dbinom{x+1}{n}_{r}.
\label{(downAll)}
\end{equation}
\end{proof}

\begin{corollary}
The number of Dyck paths to $\left( 2n,0\right) $ avoiding $r$ down steps is 
\begin{equation}
s_{n}(n)=\dfrac{1}{n+1}\dbinom{n+1}{n}_{r}.  \label{(down)}
\end{equation}
\end{corollary}

\section{Ballot paths without the pattern $\uparrow ^{r}$}

\begin{definition}
$t_{n}(m;r)=t_{n}(m)$ is the number of $\left\{ \uparrow ,\rightarrow
\right\} $ paths staying weakly above the line $y=x$ from $(0,0)$ to $(n,m)$
avoiding a sequence of $r>0$ consecutive $\uparrow $ steps. We set $%
t_{n}(m)=0$ if $n<0$ or $m+1=n>0$.
\end{definition}

This time we do not immediately have a polynomial sequence, as the table
below shows. The path $N^{r-1}\left( EN^{r-1}\right) ^{k}$ to $\left(
k,\left( r-1\right) \left( k+1\right) \right) $ is the only admissible path
reaching the point $\left( k,\left( r-1\right) \left( k+1\right) \right) $
(all others would have $r$ or more $N$-steps). Hence $t_{n-1}\left( \left(
r-1\right) n\right) =1$ for all $n\geq 1$, and $t_{n-1}\left( m\right) =0$
for $m>\left( r-1\right) n$. The only other 1's in the table occur in column
0, $t_{0}\left( m\right) =1$ for $m=0,\dots ,r-1$, and $0$ for all other
values of $m$.

The table contains a strip weakly above the diagonal $y=x$ where%
\begin{equation}
t_{n}\left( m\right) =t_{n}(m-1)+t_{n-1}\left( m\right)  \label{(RecStart)}
\end{equation}%
This happens for $0<n\leq m<n+r$ because paths in this strip cannot have $r$
consecutive vertical steps. All paths that reach a point $\left( n,m\right) $
for $m\geq n+r$ and violate the condition of not containing $N^{r}$ must
have this pattern \textbf{exactly} at the end of the path, which means that
they end in the pattern $EN^{r}$. Hence for $m\geq n+r$ we get the recurrence%
\begin{equation}
t_{n}\left( m\right) =t_{n}(m-1)+t_{n-1}\left( m\right) -t_{n-1}\left(
m-r\right)  \label{(RecEnd)}
\end{equation}%
We assume that $t_{n}\left( m\right) =0$ for all $m<n$ (also for $n=0$).

We can find a recursion that holds for all $m\geq n$ as follows: For $n\geq
1 $ we always have $t_{n}\left( n\right) =t_{n-1}\left( n\right) $, because $%
t_{n}\left( n-1\right) =0$. From (\ref{(RecStart)}) follows by induction
(inside the exceptional strip) that $t_{n}\left( m\right)
=\sum_{i=n}^{m}t_{n-1}\left( i\right) $ for all $n\leq m<n+r$. For the
values of $m$ on the boundary of the strip we have

\begin{equation*}
t_{n}\left( n+r\right) =\sum_{i=n}^{n+r}t_{n-1}\left( i\right)
-t_{n-1}\left( n\right) =\sum_{i=n+1}^{n+r}t_{n-1}\left( i\right)
\end{equation*}%
from (\ref{(RecStart)}) and (\ref{(RecEnd)}), and after that by induction
using (\ref{(RecEnd)}), \newline
$t_{n}\left( m\right) =\sum_{i=m+1-r}^{m}t_{n-1}\left( i\right) $ for all $%
m\geq n+r$. We can write both recursions together as

\begin{equation}
t_{n}\left( m\right) =\sum_{i=\max \left\{ n,m+1-r\right\}
}^{m}t_{n-1}\left( i\right)  \label{(RecSum)}
\end{equation}%
for all $m\geq n$. We can avoid the difficulty with the lower bound in the
summation by setting $t_{n}\left( m\right) =0$ for all $m\leq n$. Call the
modified numbers $t_{n}^{\prime }\left( m\right) $.The new table follows the
recursion

\begin{equation*}
t_{n}^{\prime }\left( m\right) =\sum_{i=m+1-r}^{m}t_{n-1}^{\prime }\left(
i\right)
\end{equation*}%
for all $m>n>0$. The `lost' value $t_{n}\left( n\right) $ can be easily
recovered, because $t_{n}\left( n\right) =t_{n-1}\left( n\right)
=t_{n-1}^{\prime }(n)$.

\begin{equation*}
\begin{tabular}{l}
$%
\begin{tabular}{r||rrrrrrrrr}
$m$ & 0 & 0 & 1 & 19 & 112 & 397 & 1027 & 1966 & 2905 \\ 
8 & 0 & 0 & 3 & 28 & 116 & 321 & 630 & 939 & (939) \\ 
7 & 0 & 0 & 6 & 33 & 101 & 205 & 309 & (309) & 0 \\ 
6 & 0 & 1 & 9 & 32 & 68 & 104 & (104) & 0 & 0 \\ 
5 & 0 & 2 & 10 & 23 & 36 & (36) & 0 & 0 & 0 \\ 
4 & 0 & 3 & 8 & 13 & (13) & 0 & 0 & 0 & 0 \\ 
3 & 1 & 3 & 5 & (5) & 0 & 0 & 0 & 0 & 0 \\ 
2 & 1 & 2 & (2) & 0 & 0 & 0 & 0 & 0 & 0 \\ 
1 & 1 & (1) & 0 & 0 & 0 & 0 & 0 & 0 & 0 \\ 
0 & (1) & 0 & 0 & \multicolumn{6}{r}{The values in parentheses are $0$ in $%
t_{n}^{\prime }\left( m\right) $} \\ \hline\hline
$n:$ & 0 & 1 & 2 & 3 & 4 & 5 & 6 & 7 & 8%
\end{tabular}%
$ \\ 
\multicolumn{1}{c}{The restricted ballot path counts $t_{n}\left( m\right) $
($r=4$).}%
\end{tabular}%
\end{equation*}%
In order to show the polynomial structure in the above table, we transform
it into the table below by a 90$^{\circ }$ counterclockwise turn, and
shifting the top 1's flush against the $y$-axis. In formulas, we define $%
p_{n}(m)=t_{m-1}^{\prime }((r-1)m-n)$ for $m\left( r-1\right) \geq n\geq 0$
(or $t_{n}^{\prime }\left( m\right) =p_{\left( r-1\right) \left( n+1\right)
-m}\left( n+1\right) $). The recursion $t_{n}^{\prime }\left( m\right)
=\sum_{i=m+1-r}^{m}t_{n-1}^{\prime }\left( i\right) $ `along the previous
column' becomes now a recursion $p_{n}(m)=\sum_{i=0}^{r-1}p_{n-j}(m-1)$
`along the previous row'. More precisely, for $\left( r-1\right) m-n>m-1\geq
1$ , i.e. , $m\geq 2$ and $n\leq \left( r-2\right) m$ holds%
\begin{eqnarray}
p_{n}(m) &=&t_{m-1}^{\prime }((r-1)m-n)=\sum_{i=(r-1)\left( m-1\right)
-n}^{\left( r-1\right) m-n}t_{m-2}^{\prime }\left( i\right)  \notag \\
&=&\sum_{i=(r-1)\left( m-1\right) -n}^{\left( r-1\right) m-n}p_{\left(
r-1\right) \left( m-1\right) -i}(m-1)=\sum_{i=0}^{r-1}p_{n-j}(m-1)
\label{(RecP)}
\end{eqnarray}%
The numbers $p_{n}\left( m\right) $ for $0\leq n\leq \left( r-2\right) m$
are exactly the cases where $t_{n}^{\prime }\left( m\right) $ is positive,
and the only additional numbers needed in the recursion (\ref{(RecP)}) are
the numbers \newline
$p_{\left( r-2\right) m-j}\left( m-1\right) =t_{m-2}^{\prime }((r-1)\left(
m-1\right) -\left( r-2\right) m+j)=t_{m-2}^{\prime }(m+j-r+1)=0$ for $%
j=0,\dots ,r-3$. We also add a row $p_{n}\left( 0\right) =\delta _{n,0}$ for 
$n=0,\dots ,r-2$ to the table, so that the recursion (\ref{(RecP)}) holds
for $m=1$. This part of the $p$-table, $p_{n}\left( m\right) $ for $0\leq m$
and $0\leq n\leq \left( r-2\right) \left( m+1\right) $, is shown below for $%
r=4$. Note that $p_{0}\left( m\right) =t_{m-1}((r-1)m)=1$ for all $m\geq 1$,
and also $p_{0}\left( 0\right) =1$.

\begin{equation*}
\begin{tabular}{l}
$%
\begin{tabular}{r||rrrrrrrrr}
$m$ & 1 & 8 & 36 & 119 & 315 & 699 & 1338 & 2246 & 3344 \\ 
7 & 1 & 7 & 28 & 83 & 197 & 391 & 667 & 991 & 1295 \\ 
6 & 1 & 6 & 21 & 55 & 115 & 200 & 297 & 379 & 419 \\ 
5 & 1 & 5 & 15 & 34 & 61 & 90 & 112 & 116 & \textbf{101} \\ 
4 & 1 & 4 & 10 & 19 & 28 & 33 & \textbf{32} & 23 & 13 \\ 
3 & 1 & 3 & 6 & 9 & \textbf{10} & 8 & 5 & 0 & 0 \\ \cline{9-10}
2 & 1 & 2 & \textbf{3} & 3 & 2 & 0 & 0 & \multicolumn{1}{|r}{\textit{-2}} & 
\textit{2} \\ \cline{7-8}
1 & \textbf{1} & 1 & 1 & 0 & 0 & \multicolumn{1}{|r}{\textit{-1}} & \textit{1%
} & \textit{-2} & \textit{4} \\ \cline{5-6}
0 & 1 & 0 & 0 & \multicolumn{1}{|r}{\textit{-1}} & \textit{1} & \textit{-1}
& \textit{2} & \textit{-4} & \textit{7} \\ \hline\hline
$n:$ & 0 & 1 & 2 & 3 & 4 & 5 & 6 & 7 & 8%
\end{tabular}%
$ \\ 
The rotated and shifted table $p_{n}\left( m\right) $ and its polynomial \\ 
extension below the staircase, for $r=4$. The bold numbers \\ 
occur also on the second subdiagonal in the table below.%
\end{tabular}%
\end{equation*}%
We obtained the recursion (\ref{(RecSum)}) as a discrete integral from (\ref%
{(RecEnd)}) and (\ref{(RecStart)}). We can now take differences in recursion
(\ref{(RecP)}) and get $p_{n}(m)-p_{n-1}\left( m\right)
=p_{n}(m-1)-p_{n-r}(m-1)$, or%
\begin{equation*}
p_{n}(m)-p_{n}\left( m-1\right) =p_{n-1}(m)-p_{n-r}(m-1)
\end{equation*}%
for all $m\geq 1$ and $0\leq n\leq \left( r-2\right) \left( m+1\right) $.
The column $p_{0}\left( m\right) $ can be extended as a column of ones to
all integers $m$; hence $p_{0}\left( m\right) $ can be extended to the
constant polynomial $1$. The recursion shows by induction that the $n$-th
column can be extended to a polynomial of degree $n$, and by interpolation
we can assume that we have polynomials in a real variable. The extension of $%
p_{n}\left( m\right) $ is again denoted by $p_{n}\left( m\right) $. The
above table shows some values of the polynomial expansion in cursive. The
expansion follows the same recursion, hence%
\begin{equation}
p_{n}(x)-p_{n}\left( x-1\right) =p_{n-1}(x)-p_{n-r}(x-1)  \label{(Up2)}
\end{equation}%
with initial values $p_{\left( r-2\right) m+j}\left( m\right) =0$ for $%
j=1,\dots ,r-2$ and $m\geq 0$. These conditions, together with $p_{0}\left(
0\right) =1$, determine the solution uniquely.

Recursion (\ref{(Up2)}) shows that $\left( p_{n}\left( x\right) \right) $ is
a Sheffer sequence for the same operator $B$ as the sequences $\left(
s_{n}\left( x\right) \right) $ in recursion (\ref{(OperatorRec1)}). Hence $%
p_{n}\left( x\right) $ can be written in terms of the same basis, the
Eulerian coefficients, as $s_{n}\left( x\right) $. However, the initial
values (zeroes) for $\left( p_{n}\left( x\right) \right) $ are more
difficult, because they are not on a line with positive slope. We introduce
know a Sheffer sequence $\left( q_{n}(x;\alpha )\right) _{n\geq 0}$ for the
delta operator $B$ that has roots on the parallel to the diagonal shifted by 
$\alpha +1$, $q_{n}(n-\alpha -1;\alpha )=0$, and agrees with $\left(
p_{n}\right) $ at one position left of the roots, for each $n$.

\begin{lemma}
\label{LemComposition}For the Sheffer sequence $\left( q_{n}\left( x;\alpha
\right) \right) $ for $B$ with initial values \newline
$q_{0}\left( m;\alpha \right) =1$, $q_{n}(0)=\delta _{n,0}$ for $0\leq n\leq
\alpha $ and $q_{n}(n-\alpha -1;\alpha )=0$ for $n>\alpha $, \newline
holds%
\begin{equation*}
q_{n+\alpha }\left( n;\alpha \right) =p_{\left( r-2\right) n-\alpha }\left(
n\right)
\end{equation*}%
for all $n\geq \left\lceil \alpha /\left( r-2\right) \right\rceil $.
\end{lemma}

We will proof this Lemma in Subsection \ref{SubSProof}. 
\begin{equation*}
\begin{tabular}{l}
$%
\begin{tabular}{r||rrrrrrrrr}
$m$ & 1 & 6 & 21 & 56 & 120 & 214 & 320 & 386 & \textbf{321} \\ 
5 & 1 & 5 & 15 & 35 & 65 & 99 & 121 & \textbf{101} & 0 \\ 
4 & 1 & 4 & 10 & 20 & 31 & 38 & \textbf{32} & 0 & -70 \\ 
3 & 1 & 3 & 6 & 10 & 12 & \textbf{10} & 0 & -22 & -58 \\ 
2 & 1 & 2 & 3 & 4 & \textbf{3} & 0 & -7 & -18 & -33 \\ 
1 & 1 & 1 & 1 & \textbf{1} & 0 & -2 & -6 & -10 & -15 \\ 
0 & 1 & 0 & 0 & 0 & 0 & -2 & -4 & -4 & -5 \\ \hline\hline
$n:$ & 0 & 1 & 2 & 3 & 4 & 5 & 6 & 7 & 8%
\end{tabular}%
$ \\ 
\multicolumn{1}{c}{The polynomials $q_{n}\left( m,2\right) $ for $r=4$}%
\end{tabular}%
\end{equation*}

The sequence $\left( q_{n}\right) $ agrees with the Euler coefficients $%
b_{n}\left( x\right) =\dbinom{x}{n}_{r}$ for the first degrees $n=0,\dots
,\alpha $. It follows from the Binomial Theorem for Sheffer sequences that%
\begin{equation}
q_{n}\left( x;\alpha \right) =\sum_{i=0}^{\alpha }\dbinom{i-\alpha -1}{i}_{r}%
\frac{x+\alpha +1-n}{x+\alpha +1-i}\dbinom{x+\alpha +1-i}{n-i}_{r}.
\label{(q)}
\end{equation}

\begin{corollary}
The number of ballot paths avoiding $r$ $\uparrow $-steps equals 
\begin{equation*}
t_{n}\left( m;r\right) =\sum_{i=0}^{m-n-1}\frac{1}{m+1-i}\dbinom{i-m+n}{i}%
_{r}\dbinom{m+1-i}{m-i}_{r}
\end{equation*}%
for $m>n\geq 0$. Furthermore, $t_{n}\left( n\right) =t_{n-1}\left( n\right) =%
\dfrac{1}{n+1}\dbinom{n+1}{n}_{r}$ for all $n>0$.
\end{corollary}

\begin{proof}[Proof:]
$t_{n}\left( m;r\right) =p_{\left( r-1\right) \left( n+1\right) -m}\left(
n+1\right) =q_{m}\left( n+1;m-1-n\right) $.
\end{proof}

The Corollary shows that the number of Dyck paths to $\left( 2n,0\right) $
avoiding $r$ up steps, $\dfrac{1}{n+1}\dbinom{n+1}{n}_{r}$, equals the
number of Dyck paths to $\left( 2n,0\right) $ avoiding $r$ down steps (see
formula (\ref{(down)})).

\subsection{A Conjecture for the Case $r=4$.\label{Conjecture}}

A Motzkin path can take horizontal unit steps in addition to the up and down
steps of a Dyck path. Suppose a Motzkin path is \textquotedblleft
peakless\textquotedblright , i.e., the pattern $uu$ and $ud$ does not occur
in the path. Denote the number of peakless Motzkin paths to $\left(
n,0\right) $ by $M^{\prime }\left( n\right) $. Starting at $n=0$ we get the
following sequence, $1,1,1,2,4,7,13,26,52,104,212,438,910,\dots $ for $%
M^{\prime }\left( n\right) $ (see http://www.research.att.com/\symbol{126}%
njas/sequences/A023431). It is easy to show that $M^{\prime }\left( n\right)
=\sum_{i=0}^{\left\lfloor n/3\right\rfloor }\dbinom{n-i}{2i}\frac{1}{i+1}%
\dbinom{2i}{i}$.

For the case $r=4$ we conjecture that $p_{n}(0)=\left( -1\right)
^{n}M^{\prime }\left( n-3\right) $ for all $n\geq 3$. That would imply $%
\sum_{n\geq 0}p_{n}\left( 0\right) t^{n}=\left( 3+t-\sqrt{\left( 1+t\right)
^{2}+4t^{3}}\right) /2$, and therefore%
\begin{equation*}
\sum_{n\geq 0}p_{n}\left( x\right) t^{n}=\frac{3+t-\sqrt{\left( 1+t\right)
^{2}+4t^{3}}}{2}\left( \frac{1-t^{4}}{1-t}\right) ^{x}.
\end{equation*}%
For example, the coefficient of $t^{7}$ in this power series equals \newline
$\left( x-3\right) \left(
x^{6}+24x^{5}+247x^{4}+426x^{3}-38x^{2}-2340x+6720\right) /7!$, \newline
which in turn equals $p_{7}\left( x\right) $, as can be checked using the
table for $p_{7}\left( m\right) $.

\subsection{Proof of Lemma \protect\ref{LemComposition}\label{SubSProof}}

Because of recursion (\ref{(RecP)}) we obtain the operator identity%
\begin{equation*}
I=E^{-1}\left( B^{0}+B^{1}+\dots +B^{r-1}\right)
\end{equation*}%
which holds for $\left( p_{n}\left( x\right) \right) $ and $\left(
q_{n}\left( x;\alpha \right) \right) $, and shows that both polynomials
enumerate lattice paths with steps $\left\langle 0,1\right\rangle
,\left\langle 1,1\right\rangle ,\,\left\langle 2,1\right\rangle ,\dots
,\left\langle r-1,1\right\rangle $ (above the respective boundaries). The
number of such paths reaching $\left( n,m\right) $ can also be seen as
compositions of $m$ into $n$ terms taken from $\left\{ 0,1,\dots
,r-1\right\} $. In the case of $p_{n}\left( m\right) $ the terms $%
a_{1},\dots ,a_{n}$ also have to respect the boundary, which means that $%
\sum\limits_{i=1}^{k}a_{i}\leq \left( r-2\right) k$ for all $k=1,\dots n$.
For $q_{n}\left( m;\alpha \right) $ we get for the same reason that $%
\sum\limits_{i=1}^{k}b_{i}\leq k+\alpha -1$. Such restricted compositions
have the following nice property.

\begin{lemma}
Let $c\in \mathbb{N}_{1}$, $\alpha \in \mathbb{N}_{0}$, and $n$ be a natural
number such that $n-\alpha \geq 0$. Let $P_{n}^{\alpha }$ be the number of
compositions of $cn-\alpha $ into $n$ parts from $[0,c+1]$ such that $%
a_{1}+a_{2}+\ldots +a_{n}=cn-\alpha $, and $\sum\limits_{i=1}^{k}a_{i}\leq
ck $ for all $k=1,\dots ,n-1$. Let $Q_{n}^{\alpha }$ be the number of
compositions of $n+\alpha $ into $n$ parts from $[0,c+1]$ such that $%
b_{1}+b_{2}+\ldots +b_{n}=n+\alpha $, and $\sum\limits_{i=1}^{k}b_{i}\leq
k+\alpha $ for all $k=1,\dots ,n-1$. Then $P_{n}^{\alpha }=Q_{n}^{\alpha }$.
\end{lemma}

\begin{proof}[Proof:]
Suppose, $b_{1}+b_{2}+\ldots +b_{n}=n+\alpha $, and $\sum%
\limits_{i=1}^{k}b_{i}\leq k+\alpha $. Define $a_{i}=c+1-b_{n+1-i}$ for $%
i=1,\ldots ,n$. Note that $a_{i}\in \lbrack 0,c+1]$, and $n-k+\alpha \geq
\sum\limits_{i=1}^{n-k}b_{i}=n+\alpha -\sum\limits_{i=1}^{k}b_{n+1-i}$. Hence%
\begin{equation*}
\sum\limits_{i=1}^{k}a_{i}=(c+1)k-\sum\limits_{i=1}^{k}b_{n+1-i}\leq
(c+1)k-n+(n-k)=ck
\end{equation*}%
and%
\begin{equation*}
\sum\limits_{i=1}^{n}a_{i}=(c+1)n-\sum\limits_{i=1}^{n}b_{n+1-i}=cn-\alpha .
\end{equation*}
\end{proof}

We apply this Lemma with $c=r-2$ to obtain $p_{\left( r-2\right) n-\alpha
}\left( n\right) =q_{n+\alpha }\left( n\right) $.

\subsection{Abelization}

Let $\left( b_{n}\right) $ be the basic sequence for some arbitrary delta
operator $B$, i.e., $Bb_{n}=b_{n-1}$ and $b_{n}\left( 0\right) =\delta
_{0,n} $. Every basic sequence is also a sequence of binomial type, which
means that $\sum_{n\geq 0}b_{n}\left( x\right) t^{n}=e^{x\beta \left(
t\right) }$, where $\beta \left( t\right) =t+a_{2}t^{2}+\dots $ is a formal
power series. The compositional inverse of $\beta \left( t\right) $ is the
power series that represents $B$,%
\begin{equation*}
B=\beta ^{-1}\left( D\right) =D+b_{2}D^{2}+\dots
\end{equation*}%
where $D=\partial /\partial x$ is the $x$-derivative. The Abelization of $%
\left( b_{n}\right) $ (by $a\in \mathbb{R}$) is the basic sequence $\left( 
\frac{x}{x+an}b_{n}\left( x+an\right) \right) _{n\geq 0}$ for the delta
operator $E^{-a}B$ (see \cite{ROK}). Note that with any Sheffer sequence $%
\left( s_{n}\right) $ for $B$ the sequence $\left( s_{n}\left( x+c-an\right)
\right) $ is a Sheffer sequence for $E^{a}B$. Hence \newline
$\left( \left( \frac{x+c-an}{x+c}b_{n}\left( x+c\right) \right) \right)
_{n\geq 0}$ is a Sheffer sequence for $E^{a}E^{-a}B=B$ again. Choosing $%
c=a=1 $ shows (\ref{(downAll)}).

Sheffer sequences and the basic sequence for the same delta operator are
connected by the Binomial Theorem for Sheffer sequences,%
\begin{equation*}
s_{n}\left( y+x\right) =\sum_{i=0}^{n}s_{i}\left( y\right) b_{n-i}\left(
x\right) .
\end{equation*}%
Applying this Theorem to $\left( \frac{x}{x+an}b_{n}\left( x+an\right)
\right) $ and $\left( b_{n}\left( x+c+an\right) \right) $ shows that%
\begin{equation*}
b_{n}\left( y+x+c+an\right) =\sum_{i=0}^{n}b_{i}\left( y+c+ai\right) \frac{x%
}{x+a\left( n-i\right) }b_{n-i}\left( x+a\left( n-i\right) \right) .
\end{equation*}%
Choosing $x$ as $x+\alpha +1-an$, $a=1$, $c=0$, and $y=-\alpha -1$ gives $%
b_{n}\left( x\right) =\sum_{i=0}^{n}b_{i}\left( i-\alpha -1\right) \frac{%
x+\alpha +1-n}{x+\alpha +1-i}b_{n-i}\left( x+\alpha +1-i\right) $.This is
not quite what we have in (\ref{(q)}); there the summation stops at $\alpha $%
. This effect in (\ref{(q)}) is due to the `initial values' $b_{i}\left(
i-\alpha -1\right) $ which are $0$ for $i>\alpha $.

The generating function of a Sheffer sequence $\left( s_{n}\right) $ for $B$
is of the form $\phi \left( t\right) e^{x\beta \left( t\right) }$, where $%
\phi \left( t\right) =\sum_{n\geq 0}s_{n}\left( 0\right) t^{n}$. If $%
s_{n}\left( x\right) =\frac{x+c-an}{x+c}b_{n}\left( x+c\right) =b_{n}\left(
x+c\right) -\frac{an}{x+c}b_{n}\left( x+c\right) $ then 
\begin{equation*}
\sum_{n\geq 1}\frac{n}{x+c}b_{n}\left( x+c\right) t^{n}=\frac{t}{x+c}\frac{%
\partial }{\partial t}e^{\left( x+c\right) \beta \left( t\right) }=t\beta
^{\prime }\left( t\right) e^{\left( x+c\right) \beta \left( t\right) }
\end{equation*}%
and%
\begin{equation*}
\sum_{n\geq 0}s_{n}\left( x\right) t^{n}=e^{\left( x+c\right) \beta \left(
t\right) }\left( 1-at\beta ^{\prime }\left( t\right) \right)
\end{equation*}%
If $c=a=1$ and $e^{\beta \left( t\right) }=\left( 1+t+\dots +t^{r-1}\right) $%
, then $\beta ^{\prime }\left( t\right) =\frac{1-\left( r-tr+t\right) t^{r-1}%
}{\left( 1-t^{r}\right) \left( 1-t\right) }$, and therefore%
\begin{equation}
\sum_{n\geq 0}s_{n}\left( m\right) t^{n}=\frac{\left( 1-t^{r}\right) ^{m}}{%
\left( 1-t\right) ^{m+2}}\left( rt^{r}\left( 1-t\right) +\left(
1-t^{r}\right) \left( 1-2t\right) \right) ,  \label{(GenFuncDown)}
\end{equation}%
the generating function of the number of ballot paths to $\left( n,m\right) $%
, avoiding $\rightarrow ^{r}$.

\section{Euler Coefficients\label{Euler}}

The coefficients of the polynomial%
\begin{equation*}
\left( 1+t+t^{2}+\cdots +t^{r-1}\right) ^{n}
\end{equation*}%
were considered by Euler in \cite{Euler}, where he gives the following
recurrence:%
\begin{equation*}
\dbinom{n}{k}_{r+1}=\sum\limits_{i=0}^{k/2}\dbinom{n}{k-i}\dbinom{k-i}{i}_{r}
\end{equation*}

To calculate the Euler coefficients in terms of only binomial coefficients,
we rewrite the polynomial as follows:%
\begin{eqnarray*}
\left( 1+t+t^{2}+\cdots +t^{r-1}\right) ^{n} &=&\left( \dfrac{1-t^{r}}{1-t}%
\right) ^{n} \\
&=&\sum\limits_{i=0}^{n}\dbinom{n}{i}t^{ri}(-1)^{i}\sum\limits_{j\geq 0}%
\dbinom{n+j-1}{j}t^{j}.
\end{eqnarray*}

Thus we have proven%
\begin{equation}
\dbinom{n}{k}_{r}=\sum\limits_{i=0}^{\left\lfloor k/r\right\rfloor }(-1)^{i}%
\dbinom{n}{i}\dbinom{n+k-ri-1}{k-ri}.  \label{(EulerBinomial)}
\end{equation}

Note that this identity implies $\lim_{r\rightarrow \infty }\dbinom{n}{k}%
_{r}=\dbinom{n+k-1}{k}$. Combinatorially, these are \textit{all} $\left\{
\uparrow ,\rightarrow \right\} $ paths avoiding $\rightarrow ^{r}$. This
problem occurs in Wilf's \emph{generatingfunctionology} \cite{Wilf}, Section
4.12. Note that identity (\ref{(EulerBinomial)}) implies $\lim_{r\rightarrow
\infty }\dbinom{n}{k}_{r}=\dbinom{n+k-1}{k}$. 
\begin{equation*}
\begin{tabular}{l}
$%
\begin{tabular}{l||lllllllll}
$x$ & 1 & 8 & 36 & 120 & 322 & 728 & 1428 & 2472 & 3823 \\ 
7 & 1 & 7 & 28 & 84 & 203 & 413 & 728 & 1128 & 1554 \\ 
6 & 1 & 6 & 21 & 56 & 120 & 216 & 336 & 456 & 546 \\ 
5 & 1 & 5 & 15 & 35 & 65 & 101 & 135 & 155 & 155 \\ 
4 & 1 & 4 & 10 & 20 & 31 & 40 & 31 & 20 & 10 \\ 
3 & 1 & 3 & 6 & 10 & 12 & 12 & 10 & 6 & 3 \\ 
2 & 1 & 2 & 3 & 4 & 3 & 2 & 1 & 0 & 0 \\ 
1 & 1 & 1 & 1 & 1 & 0 & 0 & 0 & 0 & 0 \\ 
0 & 1 & 0 & 0 & 0 & 0 & 0 & 0 & 0 & 0 \\ \hline\hline
$n:$ & 0 & 1 & 2 & 3 & 4 & 5 & 6 & 7 & 8%
\end{tabular}%
$ \\ 
\multicolumn{1}{c}{A table of Euler coefficients $\dbinom{x}{n}_{4}$ for $%
r=4 $}%
\end{tabular}%
\end{equation*}

We now show some properties about Euler Coefficients similar to the basic
properties of binomial coefficients.

\begin{enumerate}
\item For binomial coefficients, this property is usually called Pascal's
Identity:%
\begin{equation*}
\dbinom{n}{k}_{r}=\sum\limits_{i=0}^{r-1}\dbinom{n-1}{k-i}_{r}
\end{equation*}

\begin{proof}
\begin{eqnarray}
(1+t+\cdots +t^{r-1})^{n} &=&(1+t+\cdots +t^{r-1})^{n-1}(1+t+\cdots +t^{r-1})
\notag \\
&=&\sum\limits_{i=0}^{r-1}t^{i}(1+t+\cdots +t^{r-1})^{n-1}  \notag
\end{eqnarray}%
so 
\begin{equation}
\dbinom{n}{k}_{r}=\sum\limits_{i=0}^{r-1}[t^{k-i}](1+t+\cdots
+t^{r-1})^{n-1}=\sum\limits_{i=0}^{r-1}\dbinom{n-1}{k-i}_{r}  \notag
\end{equation}
\end{proof}

\item The table of Euler Coefficients is symmetric similar to Pascal's
Triangle:%
\begin{equation*}
\dbinom{n}{k}_{r}=\dbinom{n}{n(r-1)-k}_{r}
\end{equation*}

\begin{proof}
We proceed by induction on $n$, fixing $r\geq 2$. For $n=2$ we have the well
known symmetry for binomial coefficients. Suppose true for some $l>2$. From
the above recurrence we have%
\begin{eqnarray}
\dbinom{l+1}{k}_{r} &=&\sum\limits_{i=0}^{r-1}\dbinom{l}{k-i}%
_{r}=\sum\limits_{i=0}^{r-1}\dbinom{l}{l(r-1)-(k-i)}_{r}  \notag \\
&=&\sum\limits_{i=0}^{r-1}\dbinom{l}{(l+1)(r-1)-k-i}_{r}=\dbinom{l+1}{%
(l+1)(r-1)-k}_{r}  \notag
\end{eqnarray}%
and the induction follows.
\end{proof}

\item This property is similar to Vandermondt Convolution for binomial
coefficients:%
\begin{equation*}
\dbinom{n+m}{k}_{r}=\sum\limits_{i=0}^{k}\dbinom{n}{i}_{r}\dbinom{m}{k-i}%
_{r}.
\end{equation*}%
It follows because the Euler coefficients are of binomial type \cite{ROK}.

\item Here we have an identity that is trivial for binomial coefficients,
i.e. $r=2$, and gives and identity for the Catalan numbers as $r\rightarrow
\infty $.%
\begin{equation*}
\dfrac{1}{n+1}\dbinom{n+1}{n}_{r}=\dbinom{n}{n}_{r}-\sum\limits_{i=1}^{r-2}i%
\dbinom{n}{n-i-1}_{r}
\end{equation*}

\begin{proof}
Let $s_{n}(x)=\dfrac{x-n+1}{x+1}\dbinom{x+1}{n}_{r}$, as in (\ref{(downAll)}%
). The binomial theorem for Sheffer sequences states that $%
s_{n}(x+y)=\sum\limits_{i=0}^{n}s_{i}(y)b_{n-i}(x)$. Let $x=n$, $y=0$ and
noting that $s_{n}(0)=(1-n)\dbinom{1}{n}_{r}=1-n$ for $0<n<r$ and $0$
otherwise, we have $s_{n}(n)=\sum\limits_{i=0}^{n}s_{i}(0)b_{n-i}(n)$, hence%
\begin{eqnarray}
\dfrac{1}{n+1}\dbinom{n+1}{n}_{r} &=&\sum\limits_{i=0}^{r-1}(1-i)\dbinom{n}{%
n-i}_{r}  \notag \\
&=&\dbinom{n}{n}_{r}-\sum\limits_{i=1}^{r-2}i\dbinom{n}{n-i-1}_{r}.  \notag
\end{eqnarray}%
We have already noted that $\dbinom{n}{k}_{r}\rightarrow \dbinom{n+k-1}{k}$
as $r\rightarrow \infty $, so%
\begin{equation*}
\lim_{r\rightarrow \infty }\dfrac{1}{n+1}\dbinom{n+1}{n}_{r}=\dbinom{2n-1}{n}%
-\sum\limits_{i=1}^{n-1}i\dbinom{2n-i-2}{n-1}=C_{n}.
\end{equation*}
\end{proof}
\end{enumerate}


\begin{thebibliography}{9}
\bibitem{Deutsch} E. Deutsch, Dyck path enumeration, Discrete Math. \textbf{%
204} (1999) 167 -- 202.

\bibitem{Euler} L. Euler, De evolutione potestatis polynomialis cuiuscunque $%
(1+x+x^{2}+x^{3}+x^{4}+etc.)^{n}$, \emph{Nova Acta Academiae Scientarum
Imperialis Petropolitinae} \textbf{12} (1801) 47 -- 57.

\bibitem{Nied03} H. Niederhausen, Rota's umbral calculus and recursions, 
\emph{Algebra Univers}. \textbf{49} (2003) 435 -- 457.

\bibitem{ROK} G.-C. Rota,D. Kahaner, and A. Odlyzko, On the Foundations of
Combinatorial Theory VIII: Finite operator calculus, \emph{J. Math. Anal.
Appl.} \textbf{42} (1973) 684 -- 760.

\bibitem{Sapounakis} A. Sapounakis, L. Tasoulas, and P. Tsikouras, Counting
strings in Dyck paths, \emph{Discrete Math.}, to appear 2007.

\bibitem{Wilf} H. S. Wilf, \emph{generatingfunctionology}, third edition, A.
K. Peters, Ltd. 2006.
\end{thebibliography}
\end{document}